\documentclass[11]{article}
\usepackage{latexsym}
\usepackage{amsmath}
\usepackage{amssymb}
\usepackage{amscd}
\usepackage{graphicx}

\begin{document}

\newtheorem{Th}{Theorem}[section]
\newtheorem{Cor}{Corollary}[section]
\newtheorem{Prop}{Proposition}[section]
\newtheorem{Lem}{Lemma}[section]
\newtheorem{Def}{Definition}[section]
\newtheorem{Rem}{Remark}[section]
\newtheorem{Ex}{Example}[section]
\newtheorem{stw}{Proposition}[section]


\newcommand{\bet}{\begin{Th}}
\newcommand{\ent}{\stepcounter{Cor}
   \stepcounter{Prop}\stepcounter{Lem}\stepcounter{Def}
   \stepcounter{Rem}\stepcounter{Ex}\end{Th}}


\newcommand{\bec}{\begin{Cor}}
\newcommand{\enc}{\stepcounter{Th}
   \stepcounter{Prop}\stepcounter{Lem}\stepcounter{Def}
   \stepcounter{Rem}\stepcounter{Ex}\end{Cor}}
\newcommand{\bep}{\begin{Prop}}
\newcommand{\enp}{\stepcounter{Th}
   \stepcounter{Cor}\stepcounter{Lem}\stepcounter{Def}
   \stepcounter{Rem}\stepcounter{Ex}\end{Prop}}
\newcommand{\bel}{\begin{Lem}}
\newcommand{\enl}{\stepcounter{Th}
   \stepcounter{Cor}\stepcounter{Prop}\stepcounter{Def}
   \stepcounter{Rem}\stepcounter{Ex}\end{Lem}}
\newcommand{\bef}{\begin{Def}}
\newcommand{\enf}{\stepcounter{Th}
   \stepcounter{Cor}\stepcounter{Prop}\stepcounter{Lem}
   \stepcounter{Rem}\stepcounter{Ex}\end{Def}}
\newcommand{\ber}{\begin{Rem}}
\newcommand{\enr}{
   \stepcounter{Th}\stepcounter{Cor}\stepcounter{Prop}
   \stepcounter{Lem}\stepcounter{Def}\stepcounter{Ex}\end{Rem}}
\newcommand{\bee}{\begin{Ex}}
\newcommand{\ene}{
   \stepcounter{Th}\stepcounter{Cor}\stepcounter{Prop}
   \stepcounter{Lem}\stepcounter{Def}\stepcounter{Rem}\end{Ex}}
\newcommand{\Proof}{\noindent{\it Proof\,}:\ }

\newcommand{\EE}{\mathcal{E}}
\newcommand{\QQ}{\mathbf{Q}}
\newcommand{\R}{\mathbf{R}}
\newcommand{\C}{\mathbf{C}}
\newcommand{\ZZ}{\mathbf{Z}}
\newcommand{\KK}{\mathbf{K}}
\newcommand{\NN}{\mathbf{N}}
\newcommand{\PP}{\mathbf{P}}
\newcommand{\HH}{\mathbf{H}}
\newcommand{\uuu}{\boldsymbol{u}}
\newcommand{\xxx}{\boldsymbol{x}}
\newcommand{\aaa}{\boldsymbol{a}}
\newcommand{\bbb}{\boldsymbol{b}}
\newcommand{\AAA}{\mathbf{A}}
\newcommand{\BBB}{\mathbf{B}}
\newcommand{\ccc}{\boldsymbol{c}}
\newcommand{\iii}{\boldsymbol{i}}
\newcommand{\jjj}{\boldsymbol{j}}
\newcommand{\kkk}{\boldsymbol{k}}
\newcommand{\rrr}{\boldsymbol{r}}
\newcommand{\FFF}{\boldsymbol{F}}
\newcommand{\yyy}{\boldsymbol{y}}
\newcommand{\ppp}{\boldsymbol{p}}
\newcommand{\qqq}{\boldsymbol{q}}
\newcommand{\nnn}{\boldsymbol{n}}
\newcommand{\vvv}{\boldsymbol{v}}
\newcommand{\eee}{\boldsymbol{e}}
\newcommand{\fff}{\boldsymbol{f}}
\newcommand{\www}{\boldsymbol{w}}
\newcommand{\0}{\boldsymbol{0}}
\newcommand{\lon}{\longrightarrow}
\newcommand{\ga}{\gamma}
\newcommand{\pa}{\partial}
\newcommand{\QED}{\hfill $\Box$}
\newcommand{\id}{{\mbox {\rm id}}}
\newcommand{\OO}{{\mathcal O}}
\newcommand{\spcodim}{{\mbox {\rm sp-cod}}}
\newcommand{\diffcodim}{{\mbox {\rm diff-cod}}}
\newcommand{\sd}{{\mbox {\rm sd}}}
\newcommand{\KKer}{{\mbox {\rm Ker}}}
\newcommand{\Ker}{{\mbox {\rm Ker}}}
\newcommand{\grad}{{\mbox {\rm grad}}}
\newcommand{\ind}{{\mbox {\rm ind}}}
\newcommand{\rot}{{\mbox {\rm rot}}}
\newcommand{\diver}{{\mbox {\rm div}}}
\newcommand{\Gr}{{\mbox {\rm Gr}}}
\newcommand{\LG}{{\mbox {\rm LG}}}
\newcommand{\Diff}{{\mbox {\rm Diff}}}
\newcommand{\Symp}{{\mbox {\rm Symp}}}
\newcommand{\Ct}{{\mbox {\rm Ct}}}
\newcommand{\Uns}{{\mbox {\rm Uns}}}
\newcommand{\rank}{{\mbox {\rm rank}}}
\newcommand{\sign}{{\mbox {\rm sign}}}
\newcommand{\Spin}{{\mbox {\rm Spin}}}
\newcommand{\Sp}{{\mbox {\rm sp}}}
\newcommand{\Int}{{\mbox {\rm Int}}}
\newcommand{\Hom}{{\mbox {\rm Hom}}}
\newcommand{\Reg}{{\mbox {\rm Reg}}}
\newcommand{\codim}{{\mbox {\rm codim}}}
\newcommand{\ord}{{\mbox {\rm ord}}}
\newcommand{\Iso}{{\mbox {\rm Iso}}}
\newcommand{\corank}{{\mbox {\rm corank}}}
\def\mod{{\mbox {\rm mod.\ }}}
\newcommand{\pt}{{\mbox {\rm pt}}}
\newcommand{\cod}{{\mbox {\rm cod}}}
\newcommand{\enP}{\hfill $\Box$ \par\vspace{5truemm}}
\newcommand{\qed}{\hfill $\Box$ \par}
\newcommand{\spe}{\vspace{0.4truecm}}
\newcommand{\Tan}{{\mbox {\rm Tan}}}

\newcommand{\dint}[2]{{\displaystyle\int}_{{\hspace{-1.9truemm}}{#1}}^{#2}}


\title{
Legendre singularities of sub-Riemannian geodesics
}

\author{Goo Ishikawa and Yumiko Kitagawa}

\date{ }

\maketitle

\renewcommand{\thefootnote}{\fnsymbol{footnote}}



\begin{abstract}
\noindent
Let $M$ be a surface with a Riemannian metric and $UM$ the unit tangent bundle over $M$ with 
the canonical contact sub-Riemannian structure $D \subset T(UM)$. 
In this paper, the complete local classification of singularities, under the Legendre projection $UM \to M$, is given for sub-Riemannian geodesics of $(UM, D)$. 
Legendre singularities of sub-Riemannian geodesics are classified completely also for another Legendre projection from $UM$ to the space of Riemannian geodesics on $M$. The duality on Legendre singularities is observed related to the pendulum motion. 

\end{abstract}


\section{Introduction}

Let $M$ be a $C^\infty$ surface with a Riemannian metric $g$. Then the unit tangent bundle $UM$ over 
$M$ has the canonical contact structure $D \subset T(UM)$. Moreover $D$ has a sub-Riemannian structure 
induced from the Riemannian metric on $M$. A {\it sub-Riemannian geodesic} of $D$ (or a 
{\it $D$-geodesic}) is a curve on $UM$ 
which is tangent to $D$ and is a local minimizer of the sub-Riemannian or  
Carnot-Carath{\' e}odory arc length for the metric on $D$ (\cite{Montgomery}). 
Any $D$-geodesic on $UM$ is known to be an immersion if it is not a constant map. However 
the projection $\pi : UM \to M$, which is a Legendre projection,  restricted to a $D$-geodesic on $UM$ 
may have singularities, which are called the {\it Legendre singularities}.

In this paper we study Legendre singularities of $D$-geodesics on $(UM, D)$ and 
give the local classification result which determines the Legendre singularities of $D$-geodesics completely. 

The unit tangent bundle $UM$ has the geodesic flow for the metric $g$ on $M$ and is foliated by the horizontal lifts of Riemannian geodesics on $M$ to $UM$ 
for the projection $\pi : UM \to M$. 
Each leaf is a Legendre curve for the contact structure $D$ and then 
we have another Legendre projection $\pi'$, at least locally, from $UM$ to the leaf space, i.e. {\it the space of Riemannian geodesics}. 

We determine Legendre singularities of $D$-geodesics on $(UM, D)$ also for the projection 
$\pi'$ completely in this paper. 

\bet
\label{Main-Theorem} 
Let $\Gamma : (\R, t_0) \to UM$ be any germ of $D$-geodesic. Then the composite mapping 
diagram 
$
(\Gamma, \pi) : 
(\R, t_0) \xrightarrow{\Gamma} (UM, \Gamma(t_0)) \xrightarrow{\pi} (M, \pi(\Gamma(t_0))
$ 
is Legendre equivalent to one of following normal forms: 

{\rm (i)} $(c_1, \Pi), c_1 : (\R, 0) \to (\R^3, 0), c_1(t) = (0, 0, 0)$, 

{\rm (ii)} $(c_2, \Pi), c_2 : (\R, 0) \to (\R^3, 0), c_2(t) = (0, 0, t)$, 

{\rm (iii)} $(c_3, \Pi), c_3 : (\R, 0) \to (\R^3, 0), c_3(t) = (t, 0, 0)$, 

{\rm (iv)} $(c_4, \Pi), c_4 : (\R, 0) \to (\R^3, 0), c_4(t) = (\frac{1}{2}t^2, \frac{1}{3}t^3, t)$. 
\\
Here $\R^3$ with coordinates $(x, y, p)$ has the canonical contact structure defined by $dy - pdx = 0$ 
and $\Pi : (\R^3, 0) \to \R^2$ is the Legendre projection defined by $\Pi(x, y, p) = (x, y)$. 

Moreover the pair of Legendre equivalence classes of $(\Gamma, \pi)$ and $(\Gamma, \pi')$ is 
given by 
\\
{\rm ((i), (i))}, {\rm ((ii), (iii))}, {\rm ((iii), (ii))}, {\rm ((iii), (iii))}, {\rm ((iii), (iv))} or {\rm ((iv), (iii))}. 
\ent

In Theorem \ref{Main-Theorem}, the case (i) means that $\Gamma$ itself is a constant curve, 
(ii) (resp. (iii)) means $\Gamma$ is an embedding to a $\pi$-fiber, (resp. $\pi'$-fiber), and (iv) means that 
$\Gamma$ has the cusp singularities by the Legendre projection $\pi$ or $\pi'$. Note that 
the projection $\pi\circ\Gamma$ (resp. $\pi'\circ\Gamma$) of any $D$-geodesic $\Gamma$ 
is a front with only cusp singularities, provided it is not a constant map. 

The transformation of a Riemannian geodesic 
to the $\pi'$-projection of its $\pi$-lift is a kind of Legendre transformation. 
For instance, 
the set of oriented geodesics of the unit sphere $S^2$ in $\R^3$ is identified to itself by taking the orthogonal cuts of $S^2$ by the orthogonal planes to unit vectors in $S^2$. 
The set of oriented geodesics on the hyperbolic space modelled in the Minkowski $3$-space 
$\R^{2,1}$ is identified to the de-Sitter space $S^{1, 1}$. 
Moreover the space of geodesics on the Euclidean plane $\R^2$ is 
identified with $S^1\times \R$ naturally in the framework of projective duality(\cite{Ishikawa}). 
Then Theorem \ref{Main-Theorem} provides the complete local classification of projections to both surfaces for 
any oriented sub-Riemannian geodesics on the unit tangent bundle in each case. 

In this paper we investigate locally such Legendre transformations and related 
\lq\lq projective duality" on surfaces along the idea in sub-Riemannian contact geometry and geometric control theory, but in classical differential geometric language. 

In \S \ref{Basic constructions from Riemannian surfaces} 
we recall basic constructions related to Riemannian surfaces, 
and in \S \ref{Around the recognition of cusps}
 we recall some facts in singularities of differentiable mappings used for the proof of 
 Theorem \ref{Main-Theorem}. 
 We prove Theorem \ref{Main-Theorem} in the flat case. 
 After a preliminary from basic differential geometry of surfaces in \ref{Geodesic parallel coordinates}, 
we show Theorem \ref{Main-Theorem} 
in the general case in \S \ref{The case of general Riemannian surfaces}. 
In the last section \S \ref{Appendix: A naive motivation}, 
we mention a native motivation of our problem treated in this paper. 

For geometric control theory and sub-Riemannian geometry, consult \cite{AS, ABB, Montgomery, Morimoto, Kitagawa}. 
The sub-Riemannian geometry on $UM$ or $U^*M$, the unit cotangent bundle in the flat case $M = \R^2$ 
has been investigated in detail, in particular, 
the problems on conjugate-loci, cut-loci and wavefronts for the sub-Riemannian geodesics were solved in \cite{MS, Sachkov}. 
Though our aim in this paper to study on Legendre duality of singularites, 
our method of construction in the present paper essentially follows these preceding works. 
Then it would be an interesting problem, for example, to study global behaviors of projections for general Riemannian surfaces. 

In this paper, all manifolds and maps are supposed to be of class $C^\infty$ unless otherwise stated.

\section{Basic constructions from Riemannian surfaces}
\label{Basic constructions from Riemannian surfaces} 

Let $M$ be an oriented $2$-dimensional Riemannian manifold with metric $g$ and 
$TM$ the tangent bundle of $M$. 
Let $UM$ be the unit tangent bundle over $M$, 
$$
UM := \{ (x, v) \in TM \mid x \in M, v \in T_xM, g(v, v) = 1\}. 
$$
The bundle $\pi : UM \to M$, $\pi(x, v) = x$, is a principal $SO(2) = U(1)$ bundle and is naturally regarded as the orthonormal frame bundle over $M$. The Levi-Civita (Riemannian) connection on $M$ gives the 
decomposition 
$$
T(UM) = H \oplus V
$$
into the vertical distribution $V$ of rank $1$ and the horizontal distribution $H$ of rank $2$. 
Since $\pi$ induces an isomorphism $\pi_* : H_{(x, v)} \to T_xM$, the bundle $H$ has the induced Riemannian 
metric and the orientation. Moreover, for each $x \in M$, 
the fiber $U_xM$ of $\pi$ over $x \in M$ is regarded the unit circle of the Euclidean plane $T_xM$, and therefore the bundle 
$V$ has the induced metric, which is written as $(d\theta)^2$ using a radian angle parameter $\theta$. 
Thus $UM$ has the induced Riemannian metric $g + d\theta^2$ 
from $H$ and $V$ so that $H \perp V$. 

Note that the parameter $\theta$ itself is determined if the base point on the circle is fixed. 
Therefore if a unit vector field on an open set $\Omega \subset M$ is provided, then 
the function $\theta : \pi^{-1}(\Omega) \to \R$ is determined, which is periodic along $\pi$-fibers with period \lq\lq $2\pi$". 

For each $(x, v) \in UM$, $T_xM$ is decomposed as $\langle v\rangle_{\R}\oplus \langle Jv\rangle_{\R}$, 
where $J$ is the $90^\circ$ rotation, 
and therefore we have the decomposition $H = K\oplus L$ induced by $\pi_*$. Note also that 
$L = K^\perp$ in $H$ and the orthogonal decomposition $T(UM) = K \oplus L \oplus V$. 

Recall that the {\it connection form} $\omega$ on $UM$ is characterized as the $SO(2)$-invariant $1$-form $\omega$ satisfying $\Ker(\omega) = H$ and $\omega(\xi) = 1$ for the unit tangent vector $\xi$ 
with positive direction along the $\pi$-fiber, i.e. the fundamental vector corresponding to $1 \in \R = {\mathfrak so}(2)$ (see \cite{ST}).

The {\it canonical bundle} $D \subset T(UM)$ is defined by 
$$
D := \{ (x, v; \xi) \in T(UM) \mid (x, v) \in UM, \xi \in T_{(x, v)}UM, \pi_*(\xi) \in \langle v \rangle_{\R} \}. 
$$
The distribution $D$ is a contact distribution on $UM$. Note that $D = K\oplus V$ 
and $K = D \cap H$. 

Note that the geodesic flow on $UM$ induced by the Riemannian metric of $M$ preserves $K, L$ and $V$ respectively and its trajectories, the horizontal lifts of Riemannian geodesics are integral curves of $K$. 

Recall that a contact structure on a manifold $W$ means a subbundle $D \subset TW$ of codimension $1$ such that, 
any local $1$-form $\alpha$ on $W$ defining $D$, satisfies that $d\alpha\vert_D$ is non-degenerate. 
Then the dimension of $W$ is odd, say, $2n+1$ for some $n$. 
An immersion $\Gamma : N \to W$ from an $n$-dimensional manifold $N$ is called a {\it Legendre immersion} if 
$d\Gamma(TN) \subset D$. A submersion $\pi : W \to M$ to an $(n+1)$-dimensional manifold $M$ is calle a {\it Legendre projection} 
if the tangent bundle of any $\pi$-fiber is contained in $D$ (\cite{Arnold, AGV}). 
In this paper we concern only on the case $n = 1$.

\bef
{\rm
A {\it pseudo-product sub-Riemannian contact structure} $D$ on a $3$-dimensional manifold $W$ is a sub-Riemannian contact 
structure $D \subset TW$ with an orthogonal decomposition $D = K\oplus V$ into subbundles $K$ and $V$ of rank $1$ respectively. 
}
\enf

Therefore $D \subset T(UM)$ is a pseudo-product sub-Riemannian contact structure on $UM$. 

Let us denote by $N := UM/K$ the local leaf space of $K$ at a point $(x, v) \in UM$ 
and $\pi' : UM \to UM/K$ the projection. 
Then we have locally the double Legendre projection
\[
M \stackrel{\pi}{\longleftarrow} UM \stackrel{\pi'}{\longrightarrow} N, 
\]
for the contact structure $D$ on $UM$. Note that $\Ker(\pi_*) = V$ and $\Ker(\pi'_*) = K$. 

For the general theory of pseudo-product structures or double Legendre projections, see \cite{IM, SY1, SY2, Takeuchi}. 

\section{Around the recognition of cusps}
\label{Around the recognition of cusps}

Recall that two composite mapping diagrams 
$(\R, t_0) \xrightarrow{\Gamma} (W, \Gamma(t_0)) \xrightarrow{\pi} (M, \pi(\Gamma(t_0))$ and 
$(\R, t'_0) \xrightarrow{\Gamma'} (W', \Gamma'(t'_0)) \xrightarrow{\pi'} (M', \pi'(\Gamma'(t'_0))$, 
where $W, W'$ are contact manifolds, 
are called {\it Legendre equivalent} if there exist diffeomorphism-germs $\sigma : (\R, t_0) \to (\R, t'_0), 
\tau : (M, \pi(\Gamma(t_0))) \to (M', \pi'(\Gamma'(t'_0))$, and a contactomorphism-germ 
$\Phi : (W, \Gamma(t_0)) \to (W', \Gamma'(t'_0))$ such that the diagram
\[
\begin{array}{ccccc}
(\R, t_0) & \stackrel{\Gamma}{\longrightarrow} & (W, \Gamma(t_0)) & \stackrel{\pi}{\longrightarrow}  & 
(M, \pi(\Gamma(t_0))) 
\\
\sigma\downarrow & & \Phi\downarrow & & \tau\downarrow
\\
(\R, t'_0) &  \stackrel{\Gamma'}{\longrightarrow} & (W', \Gamma'(t'_0)) & \stackrel{\pi'}{\longrightarrow}  & (M', \pi'(\Gamma'(t'_0)))
\end{array}
\]
is commutative (\cite{Arnold, AGV, AGLV}). Then the compositions $\pi\circ \Gamma$ and $\pi'\circ\Gamma'$ are 
right-left equivalent by diffeomorphisms $\sigma$ and $\tau$. 

A map-germ $\gamma : (\R, t_0) \to M$ to a surface is called a {\it cusp} if 
$\gamma$ 
is right-left equivalent to the standard cusp $(\R, 0) \to (\R^2, 0), t \mapsto (\frac{1}{2}t^2, \frac{1}{3}t^3)$. 

We use the following fundamental recognition lemma on cusp singularities. 

\bel
\label{recognition-of-cusp}
{\rm (\cite{USY})} 
Let $k = (x_1, x_2) : (\R, t_0) \to \R^2$ be a germ of $C^\infty$ curve on the plane. 
Suppose $k$ is not an immersion at $t_0$, i.e. $(\dot{x}_1(t_0), \dot{x}_2(t_0)) = (0, 0)$. Then 
$k$ is a cusp if and only if 
$$
\Delta := 
\left|
\begin{array}{cc}
\ddot{x}_1 & \dddot{x}_1
\\
\ddot{x}_2 & \dddot{x}_2
\end{array}
\right|(t_0) \not= 0. 
$$
\enl

\Proof
Suppose $k$ is not an immersion at $t_0$.  
Then we see, by simple direct calculations, that the condition $\Delta \not= 0$ depends only on the right-left equivalence class of $k$. 
Then we see if $k$ is cusp then $\Delta \not= 0$ for the normal form of cusp. 
Now suppose $\Delta \not= 0$. 
Then we have 
$(X_1\circ k)(T) = T^m, (X_2\circ k)(T) = T^{m+1} a(T)$ for an integer $m \geq 2$ and a $C^\infty$ function-germ $a : (\R, 0) \to \R$, 
by taking a new coordinate $T = t - t_0$ of $\R$ and a system of coordinates $(X_1, X_2)$ 
on $(\R^2, k(t_0))$ centered at $k(t_0)$. Since $\Delta \not= 0$, we have $m = 2$ and $a(0) \not= 0$. 
Then $(X_1\circ k)(T) = T^2, (X_2\circ k)(T) = T^3 a(T)$. We see there exist $C^\infty$ function-germs $b(T), c(T)$ such that 
$a(T) = b(T^2) + Tc(T^2)$. Then $(X_2\circ k)(T) = T^3 a(T^2) + T^4 c(T^2)$, using Malgrange preparation theorem 
(see \cite{BL, Martinet}). 
Set $Y_1 = X_1, Y_2 = \frac{1}{a(X_1)}\left( X_2 - X_1^2 c(X_1)\right)$. 
Then the Jacobian $\frac{\pa(Y_1, Y_2)}{\pa(X_1, X_2)} \not= 0$ at $(0, 0)$ 
and $(Y_1\circ k)(T) = T^2, (Y_2\circ k)(T) = T^3$ for the new system of coordinates $(Y_1, Y_2)$. 
After a linear transformation, we have the result. 
\QED

\ber
{\rm
It is known that any two Legendre projections $\pi : (W, z_0) \to (M, x_0)$ and $\pi' : (W', z'_0) \to (M', x'_0)$ are 
Legendre equivalent, i.e. there exist a diffeomorphism-germ $\tau : (M, x_0) \to (M', x'_0)$ and a contactomorphism-germ 
$\Phi : (W, z_0) \to (W', z'_0)$ such that $\tau\circ\pi = \pi'\circ\pi'$ (\cite{Arnold, AGV}). 
}
\enr

\bel
\label{recognition-of-cusp-2}
Let $W$ be a $3$-dimensional contact manifold, 
$\pi : W \to M$ a Legendre projection and $\Gamma : (\R, t_0) \to W$ a Legendre immersion. 
Suppose $\pi\circ\Gamma$ is not an immersion at $t_0$. 
Then we have that $\pi\circ\Gamma : (\R, t_0) \to M$ is a cusp if and only if 
the second derivative $(\pi\circ\Gamma)''(t_0) \not= 0$. 
\enl

\Proof
Assume $\pi\circ\Gamma$ is a cusp. Take a system of local coordinates $x_1, x_2, x_3$ of $W$ centered at $\Gamma(t_0)$ 
such that $x_1$ and $x_2$ are constant along each $\pi$-fibers. Then we have that $(x_1,  x_2)$ 
induces a system of local coordinates 
of $M$ centered at $\pi\circ\Gamma(t_0)$ and $\pi'$ is given by 
$(x_1, x_2, x_3) \mapsto (x_1, x_2)$. Then, by Lemma \ref{recognition-of-cusp}, 
$(x_1''(t_0), x_2''(t_0)) \not= (0, 0)$, 
for the system of local coordinates $(x_1, x_2, x_3)$. 
Conversely assume $\pi\circ\Gamma(t_0) \not= 0$. 
Then the planer curve $(x_1(\Gamma(t)), x_2(\Gamma(t))$ is singular at $t_0$ 
and a non-vanishing term of second order 
for the coordinate $T = t - t_0$. Changing the systems of local coordinates $(x_1, x_2)$ and $T$ if necessary, we have 
$x_1\circ\Gamma(T) = T^2, x_2\circ\Gamma(T) = cT^3 + e(T)$, the order of $e(T)$ at $0$ being $> 3$, for some $c \in \R$. 
Since Legendre lift of the planer curve is unique and must be an immersion which is Legendre equivalent to $\Gamma$ 
at $t_0$, we have $c \not= 0$. Then, by Lemma \ref{recognition-of-cusp} or a direct argument as in the proof of 
Lemma \ref{recognition-of-cusp}, we see that $\pi\circ\Gamma$ is a cusp, i.e., it is right-left equivalent to the normal form of the cusp. 
\QED

\section{The flat case}
\label{The flat case}

First we consider the case $M = \R^2$, the Euclidean plane with coordinates $x_1, x_2$. 
Then $UM$ has coordinates $x_1, x_2, \theta$, where $\theta$ is the radian angle coordinate for the section $\frac{\pa}{\pa x_1}$. 
We explain the general basic constructions in sub-Riemannian geometry along this simple situation. 

We set 
$$
V_1 = \cos\theta\frac{\pa}{\pa x_1} + \sin\theta\frac{\pa}{\pa x_2}, 
\quad
V_2 = \frac{\pa}{\pa \theta}, 
$$
which form an orthonormal frame of $D \subset T(UM)$. 
Let $\Gamma : [a, b] \to UM$ be an absolutely continuous or a piecewise smooth curve 
such that $\Gamma'(t) \in D$ for almost every $t \in [a, b]$. 
The sub-Riemannian or Carnot-Caratheodory arc length of $\Gamma$ is defined by
$$
L(\Gamma) = \dint{a}{b} \Vert \Gamma'(t)\Vert dt
$$
using the norm of the sub-Riemannian metric on $D$ introduced in \S \ref{Basic constructions from Riemannian surfaces}. 
It is known the length minimizing problem is equivalent to the energy minimizing problem (\cite{Montgomery}). 

We represent vectors in $D \subset T(UM)$ using the frame $V_1, V_2$ as 
$$
F(x_1, x_2, \theta; u_1, u_2) = u_1V_1 + u_2V_2 = u_1\left( \cos\theta\frac{\pa}{\pa x_1} + \sin\theta\frac{\pa}{\pa x_2}\right) + u_2\frac{\pa}{\pa\theta}.
$$
The parameters $u_1, u_2$ are regarded as {\it control parameters}. 
The {\it energy function} $E : D \to \R$ is given by using the squared norm of $F$ as
$$
E(x_1, x_2, \theta; u_1, u_2) := \frac{1}{2}(u_1^2 + u_2^2). 
$$

Now we consider the optimal control problem on $D$-integral curves of minimizing the energy $E$. 
Then the Hamiltonian function $H : D\times_{UM}T^*(UM) \to \R$ of the optimal control problem 
is given by 
$H(x, v, p) := \langle p, F(v)\rangle + cE(x, v)$. 
for some constant $c$. Here $(x, v) \in D, (x, p) \in T^*(UM)$ and $x \in M$. 
In coordinates, it is written as 
$$
H(x_1, x_2, \theta; u_1, u_2; p_1, p_2, \varphi) 
:= u_1(p_1\cos\theta + p_2\sin\theta) + u_2\varphi + \frac{1}{2}c(u_1^2 + u_2^2). 
$$

By the Pontryagin principle, any solution $(x_1(t), x_2(t), \theta(t), u_1(t), u_2(t))$ 
of the optimal control problem is obtained by the 
constrained Hamilton equation
$$
\dot{x}_1 = \frac{\pa H}{\pa p_1}, \ \dot{x}_2 = \frac{\pa H}{\pa p_2}, \ \dot{\theta} = \frac{\pa H}{\pa \varphi}, \ \ 
\dot{p}_1 = - \frac{\pa H}{\pa x_1}, \ \dot{p}_2 = - \frac{\pa H}{\pa x_2}, \ \dot{\varphi} = - \frac{\pa H}{\pa \theta}
$$
with constraint $\dfrac{\pa H}{\pa u_1} = 0, \dfrac{\pa H}{\pa u_2} = 0$ for some 
$(p_1(t), p_2(t), \varphi(t)) \not = 0, c \in \R$. 

The extremal is called {\it abnormal} if $c = 0$ and is called {\it normal} if $c \not= 0$. See \cite{IKY}. 

A curve $\Gamma : (\R, t_0) \to UM$, $\Gamma(t) = (x_1(t), x_2(t), \theta(t))$, 
is called a {\it $D$-geodesic} if the above constrained Hamiltonian equation is satisfied 
for some $p_1(t), p_2(t), \varphi(t), u_1(t), u_2(t)$ and $c \in \R$. 

In our case the condition is given by 
explicitly 
\[
\begin{aligned}
\dot{x}_1 = u_1\cos\theta, \ \dot{x}_2 = u_1\sin\theta, \ \dot{\theta} = u_2, \ \ 
\dot{p}_1 = 0, \ \dot{p}_2 = 0, \ \dot{\varphi} = u_1(p_1\sin\theta - p_2\cos\theta), 
\\
p_1\cos\theta + p_2\sin\theta + cu_1 = 0, \ \ \varphi + cu_2 = 0, \ \ (p_1, p_2, \varphi, c) \not= 0, c \in \R. 
\end{aligned}
\]

By the above condition we see that each of $p_1$ and $p_2$ is a locally constant on $t$. 
Because our distribution is a contact structure it is known that there are no non-trivial abnormal extremals. 
Here a trivial extremal means a locally constant $(x_1(t), x_2(t), \theta(t))$. 
To make sure we will check that fact in our simple situation: Suppose there exists an extremal with $c = 0$.  Then 
$p_1\cos\theta + p_2\sin\theta = 0$ and $\varphi = 0$. For any $t$ with $u_1(t) \not= 0$, 
we have $p_1\sin\theta - p_2\cos\theta = 0$. 
Then we have $p_1 = p_2 = \varphi = 0$, which leads a contradiction. Therefore $u_1(t)$ must be $0$ 
almost everywhere. Since $(p_1, p_2) \not= (0, 0)$ and it is a locally constant vector, 
we have also $(\cos\theta, \sin\theta)$ and so $\theta$ must be a locally constant, 
which implies $u_2(t) = 0$ a.e. also, which means that the extremal is trivial. 

Now suppose $c \not= 0$ and seek normal extremals. 
Then, by replacing $-\frac{1}{c}p_1, -\frac{1}{c}p_2, -\frac{1}{c}\varphi$ by $p_1, p_2, \varphi$ respectively, 
we may set $c = -1$. 
Then $u_1 = p_1\cos\theta + p_2\sin\theta, u_2 = \varphi$. 
Therefore the extremal 
$(x_1(t), x_2(t), \theta(t), p_1(t), p_2(t), \varphi(t))$ satisfies a system of ordinary differential equations 
\[
\begin{aligned}
& \dot{x}_1 = (p_1\cos\theta + p_2\sin\theta)\cos\theta, \quad
\dot{x}_2 = (p_1\cos\theta + p_2\sin\theta)\sin\theta, \quad 
\dot{\theta} = \varphi, 
\\
& \dot{p}_1 = 0, \quad \dot{p}_2 = 0, \quad \dot{\varphi} = (p_1\cos\theta + p_2\sin\theta)(p_1\sin\theta - p_2\cos\theta), 
\end{aligned}
\]
with $C^\infty$ right hand sides, and any solution is of class $C^\infty$. 
Suppose that $\Gamma$ is not an immersion at $t_0$. Then $(\dot{x}_1(t_0), \dot{x}_2(t_0), \dot{\theta}(t_0)) 
= (0, 0, 0)$. Then $p_1, p_2, \varphi$ must be all identically zero, and $\Gamma$ should be a constant map. 
Therefore any non-constant $D$-geodesic $\Gamma$ is an immersion. 
Moreover we observe that $\theta$ satisfies the second order ordinary differential equation 
\[
\begin{aligned}
\ddot{\theta} & = (p_1\cos\theta + p_2\sin\theta)(p_1\sin\theta - p_2\cos\theta)
\\
& = 
p_1^2\cos\theta\sin\theta - p_1p_2\cos^2\theta + p_1p_2\sin^2\theta - p_2^2\cos\theta\sin\theta. 
\end{aligned}
\]

Suppose the constants $p_1 = p_2 = 0$, then $u_1 = 0$ and $\dot{x_1} = \dot{x_2} = 0$, each of $x_1, x_2$ being a constant. 
Moreover $\dot{\varphi} = 0$ and $\varphi$ is a constant. 
Therefore $\dot{\theta}$ is a constant and $\theta(t) = at$ for some $a \in \R$. 
If $a = 0$, then $\Gamma$ is a constant curve. If $a \not= 0$, then 
$\Gamma$ gives a parametrization of a $\pi$-fiber over a point on $M$, and 
the $D$-geodesic $\Gamma$ can be regarded as a constant {\it directed} curve, or a frontal, on the plane 
endowed with rotating directions of constant angular velocity. 

Next suppose $(p_1, p_2) \not= (0, 0)$. 
For example, for $p_1 = 0, p_2 = 1$, then the $D$-geodesic $\Gamma(t) = (\dot{x}_1(t), \dot{x}_2(t), \theta(t))$ 
satisfies $\dot{x}_1 = \sin\theta\cos\theta, \dot{x}_2 = \sin^2\theta$ 
and $\ddot{\theta} = - \sin\theta\cos\theta = -\frac{1}{2}\sin 2\theta$. 
In general we have 
$\ddot{\theta} = - r\sin(2\theta + \rho), $
where we set $r = \frac{1}{2}(p_1^2 + p_2^2)$, $\cos \rho = -\frac{1}{2r}(p_1^2 - p_2^2)$ and 
$\sin \rho = \frac{1}{r}{p_1p_2}$. 
If we set $\Theta = 2\theta + \rho$, $\omega = \sqrt{2r}$, then we have 
$$
\ddot{\Theta} = - \omega^2\sin\Theta, 
$$
that is the non-linear equation of a {\it simple pendulum}. See also \cite{MS, Sachkov}. 
Therefore $\theta$ can be expressed by 
elliptic functions. We need only the simple behavior of $\theta$ hereafter: 
When $\dot{\theta} = 0$, then $\ddot{\theta} \not= 0$. When $\ddot{\theta} = 0$, then 
$\dot{\theta} \not= 0$.

\

Now we begin to show Theorem \ref{Main-Theorem}.

\

\noindent
{\it Proof of Theorem \ref{Main-Theorem} in the flat case.} 
Let $\Gamma : (\R, t_0) \to (UM, \Gamma(t_0))$ be a $D$-geodesic. Set $\Gamma(t) = (x_1(t), x_1(t), \theta(t))$ as above. 

{\it Part I}. \  $\pi$-Legendre classification of $\Gamma$.

We have by setting $c = -1$, 
$$
\dot{x}_1 = (p_1\cos\theta + p_2\sin\theta)\cos\theta, 
\quad 
\dot{x}_2 = (p_1\cos\theta + p_2\sin\theta)\sin\theta. 
$$

Let 
$\widetilde{\Gamma}(t) = (\Gamma(t); u_1(t), u_2(t); p_1(t), p_2(t), \varphi(t))$ 
be a corresponding extremal. 
If $\Gamma$ is a constant curve, then $(\Gamma, \pi)$ is Legendre equivalent to the case (i) of Theorem \ref{Main-Theorem}. 
If $\gamma$ is not a constant curve, but $\pi\circ\Gamma = (x_1, x_2)$ is a constant curve. Then 
$(\Gamma, \pi)$ is Legendre equivalent to (ii). 
Suppose $\pi\circ\Gamma$ is not a constant curve. 
First suppose $(\dot{x}_1(t_0), \dot{x}_2(t_0)) \not= (0, 0)$, i.e. $\pi\circ\gamma$ is an immersion-germ. 
Then $(\Gamma, \pi)$ is Legendre equivalent to (iii). 
Now suppose $(\dot{x}_1(t_0), \dot{x}_2(t_0)) = (0, 0)$. Note that, then, we have 
$p_1\cos\theta(t_0) + p_2\sin\theta(t_0) = 0$ and therefore $\ddot{\theta}(t_0) = 0$. 
Then we have 
\[
\begin{aligned}
\ddot{x}_1 & = \dot{\theta}\{ (-\sin\theta)(p_1\cos\theta + p_2\sin\theta) + \cos\theta(-p_1\sin\theta + p_2\cos\theta)\}
\ = \dot{\theta}(- p_1\sin 2\theta + p_2\cos 2\theta), 
\\
\ddot{x}_2 & = \dot{\theta}\{ \cos\theta(p_1\cos\theta + p_2\sin\theta) + \sin\theta(-p_1\sin\theta + p_2\cos\theta)\}
\ = \dot{\theta}(p_1\cos 2\theta + p_2\sin 2\theta), 
\\
\dddot{x}_1 & = \ddot{\theta}(-p_1\sin 2\theta + p_2\cos 2\theta) + \dot{\theta}^2(-2p_1\cos 2\theta - 2p_2\sin 2\theta), 
\\
\dddot{x}_2 & = \ddot{\theta}(p_1\cos 2\theta + p_2\sin 2\theta) + \dot{\theta}^2(-2p_1\sin 2\theta + 2p_2\cos 2\theta). 
\end{aligned}
\]
Therefore we have
\[
\begin{aligned}
\left|
\begin{array}{cc}
\ddot{x}_1 & \dddot{x}_1
\\
\ddot{x}_2 & \dddot{x}_2
\end{array}
\right|
& = 
\left|
\begin{array}{cc}
\dot{\theta}(- p_1\sin 2\theta + p_2\cos 2\theta) & \ddot{\theta}(-p_1\sin 2\theta + p_2\cos 2\theta) + \dot{\theta}^2(-2p_1\cos 2\theta - 2p_2\sin 2\theta)
\\
\dot{\theta}(p_1\cos 2\theta + p_2\sin 2\theta) & \ddot{\theta}(p_1\cos 2\theta + p_2\sin 2\theta) + \dot{\theta}^2(-2p_1\sin 2\theta + 2p_2\cos 2\theta)
\end{array}
\right|
\\
& = 
2\dot{\theta}^3
\left|
\begin{array}{cc}
-p_1\sin 2\theta + p_2\cos 2\theta & - p_1\cos 2\theta - p_2\sin 2\theta 
\\
p_1\cos 2\theta + p_2\sin 2\theta & -p_1\sin 2\theta + p_2\cos 2\theta
\end{array}
\right|
\\
& = 2\dot{\theta}^3\{ (-p_1\sin 2\theta + p_2\cos 2\theta)^2 + (p_1\cos 2\theta + p_2\sin 2\theta)^2\} 
= 2\dot{\theta}^3(p_1^2 + p_2^2)
\end{aligned}
\]
Thus we have $\Delta = 2\dot{\theta}(t_0)^3(p_1^2 + p_2^2)$. 
Since $\ddot{\theta}(t_0) = 0$, and
since $\theta(t)$ satisfies the above second order ordinary differential equation and is not a constant, 
we see that $\dot{\theta}(t_0) \not= 0$. 
Therefore $\Delta \not= 0$, and we see that $\pi\circ\Gamma$ is right-left equivalent to the cusp $t \mapsto 
(\frac{1}{2}t^2, \frac{1}{3}t^3)$, which has the unique Legendre lift $t \mapsto (\frac{1}{2}t^2, \frac{1}{3}t^3, t)$ 
to the standard contact manifold $\R^3$ with coordinates $(x, y, p)$ with $dy - pdx = 0$. Thus we have that 
$(\Gamma, \pi)$ is Legendre equivalent to (iv). 

{\it Part II}. \ $\pi'$-Legendre classification of $\Gamma$.

In our flat case, the projection $\pi'$ is given by 
$(x_1, x_2, \theta) \mapsto (F, E)$, where 
$$
F = - x_1 \sin\theta + x_2\cos\theta, \quad E = \theta. 
$$
Note that $F$ and $E$ are independent first integrals of the geodesic flow for the flat metric on $M = \R^2$. 
We set $f = F\circ \Gamma, e = E\circ \Gamma$. 
Then we have 
\[
\begin{aligned}
\dot{f} & = - \dot{\theta}(x_1\cos\theta + x_2\sin\theta), 
\\
\ddot{f} & = 
-\dot{\theta}(p_1\cos\theta + p_2\sin\theta) - \dot{\theta}^2(-x_1\sin\theta + x_2\cos\theta) 
- \ddot{\theta}(- x_1\cos\theta + a_2\sin\theta) 
\\
\dddot{f} 
& = 
\dot{\theta}^2(p_1\cos\theta - p_2\sin\theta) - 2\ddot{\theta}(p_1\cos\theta + p_2\sin\theta) 
+ 3\dot{\theta}\ddot{\theta}(x_1\cos\theta - x_2\sin\theta) 
\\
& \quad \quad 
+ \dot{\theta}^3(x_1\cos\theta + x_2\sin\theta) 
- \dddot{\theta}(x_1\cos\theta + x_2\sin\theta)
\end{aligned}
\]
If $\Gamma$ is a constant curve, then $(\Gamma, \pi')$ is Legendre equivalent to (i). 
If $\theta$ is a constant function, then $(\Gamma, \pi')$ is Legendre equivalent to (ii). 
If $\dot{\theta}(t_0) \not= 0$, then $\pi'\circ\Gamma$ is an immersion at $t_0$ and 
$(\Gamma, \pi')$ is Legendre equivalent to (iii). 
Suppose  $\pi'\circ\Gamma$ is not an immersion at $t_0$. Then $\dot{\theta}(t_0) = 0$. 
Then we have that
$$
\left|
\begin{array}{cc}
\ddot{f} & \dddot{f}
\\
\ddot{e} & \dddot{e}
\end{array}
\right|(t_0)
= 2\ddot{\theta}(t_0)^2(p_1\cos\theta(t_0) + p_2\sin\theta(t_0)). 
$$
If $p_1\cos\theta(t_0) + p_2\sin\theta(t_0) = 0$, then $\dot{x}_1(t_0) = \dot{x}_2(t_0) = \dot{\theta}(t_0) = \dot{p}_1(t_0) = 
\dot{p}_2(t_0) = \dot{\varphi}(t_0) = 0$, which leads that $\Gamma$ is a constant curve. 
Thus, if $\Gamma$ is not a constant curve, then we see $p_1\cos\theta(t_0) + p_2\sin\theta(t_0) \not= 0$. 
Therefore we see that $\Delta = \ddot{\theta}(t_0)^2(p_1\cos\theta(t_0) + p_2\sin\theta(t_0)) \not= 0$ whenever $\dot{\theta}(t_0) = 0$. 
Thus we have that, in this case, $(\Gamma, \pi')$ is Legendre equivalent to (iv). 

The last claim on the combination of Legendre singularities for $\pi$ and $\pi'$ is obtained just by observing the 
\lq\lq pendulum" duality on points $t = t_0$ where $\dot{\theta}(t_0) = 0, \ddot{\theta}(t_0) \not= 0$ and 
points $t = t_1$ where $\dot{\theta}(t_1) \not= 0, \ddot{\theta}(t_1) = 0$, which appears as the Legendre duality 
in our case. 
\QED

\ber
{\rm 
It is the geometry of 
the curve $\pi\circ\Gamma(t) = (x_1(t), x_2(t))$ on $\R^2$, which is the projection to $\R^2$ of a $D$-geodesic 
$\Gamma(t) = (x_1(t), x_2(t), \theta(t))$. 
We see $\pi\circ\Gamma$ is singular at $t = t_0$ when $p_1\cos\theta(t_0) + p_2\sin\theta(t_0) = 0$. 
For instance we see the curvature of the plane curve $\pi\circ\Gamma(t) = (x_1(t), x_2(t))$ is given by 
$$
\kappa(t) = \frac{\dot{\theta}(t)}{\vert p_1\cos\theta(t) + p_2\sin\theta(t)\vert }, 
$$
by simple calculations, if $\pi\circ\Gamma$ is an immersion at $t$. 
Therefore we see $\pi\circ\Gamma(t)$ has an inflection point at $t = t_0$ 
if $\dot{\theta}(t_0) = 0$. 
Moreover we have that, if $\pi\circ\Gamma$ has a cusp at $t = t_0$, then the {\it cuspidal curvature} $\kappa_c$ 
of $\pi\circ\Gamma$ at $t = t_0$ is given by 
$$
\kappa_c = 
2 \, (\sign\, \dot{\theta})\, \dfrac{\vert{\dot{\theta}}\vert^{\frac{1}{2}}}{(p_1^2 + p_2^2)^{\frac{1}{4}}}. 
$$
For the cuspidal curvature see \cite{USY}. 

Further we observe that, 
for any non-constant solution of the equation of pendulum, the points $t$ where $\dot{\theta}(t) = 0$ and $\ddot{\theta}(t) = 0$ 
appear alternately. Note that, for any $D$-geodesic $\Gamma$, 
we have an inflection point $t = t_0$ where $\theta(t_0) = 0$ and a cusp point $t = t_1$ 
on $\pi\circ\Gamma$ where $p_1\cos\theta(t_0) + p_2\sin\theta(t_0) = 0$ and so $\ddot{\theta}(t_0) = 0$, 
and $\theta(t_0) \not= 0$. 
For a sub-Riemannian geodesic $\Gamma : \R \to U\R^2$, if the variation of the angle $\theta(t)$ is small, then the inflection points, where $\dot{\theta} = 0$, and the cusp points, where $\ddot{\theta} = 0$, 
appear alternately along the non-constant projection $\pi\circ\Gamma$, which may be called a \lq\lq zigzag" curve (\cite{USY}). 
Moreover we observe the directions of cusps are all parallel. 
This will provide a severe restriction on the front curve $\pi\circ\Gamma$. 
An example of the projection of $D$-geodesic is illustrated roughly like as follows. 

\begin{center}
\includegraphics[width=10truecm, height=2truecm, clip, 
]{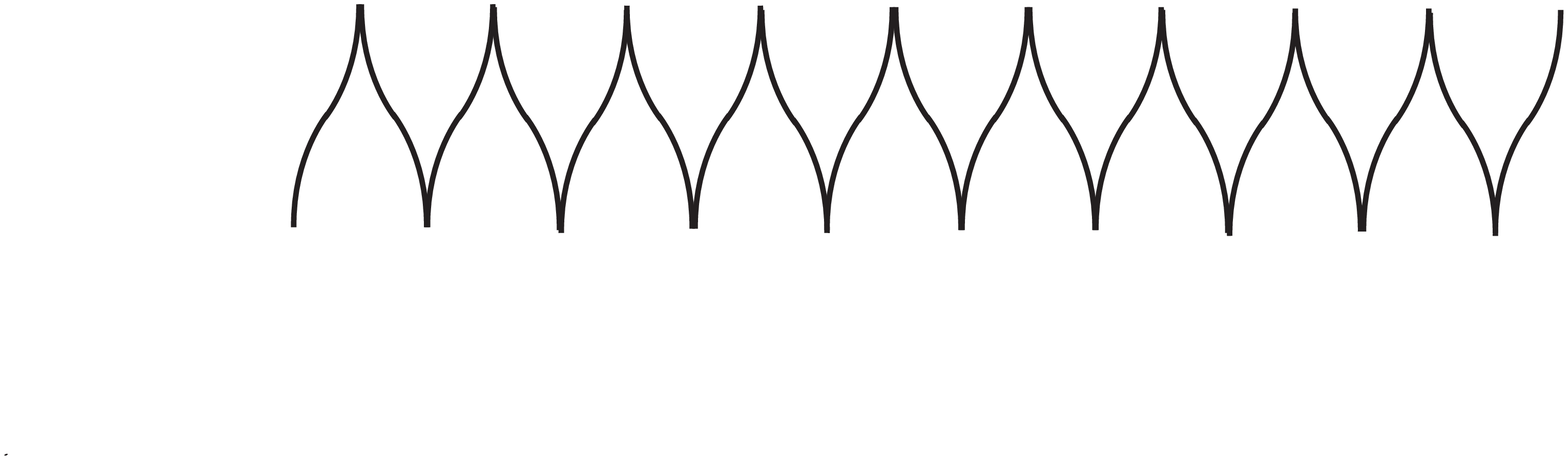} 
\end{center}
See illustrations also in \cite{ABB, MS, Sachkov}. 
}
\enr

\section{Geodesic parallel coordinates}
\label{Geodesic parallel coordinates}

To analyze the equation of sub-Riemannian geodesics on $UM$ in the case of general Riemannian surface $M$, 
it is useful to take special local coordinates on $M$, which is called {\it geodesic parallel coordinates}. 
We will recall it. 

Let $(M, g)$ be a $2$-dimensional Riemannian manifold, $p \in M$ and 
$v \in T_pM$ a unit tangent vector. 
For a system of local coordinates $(x_1, x_2)$, we set 
$g(\frac{\pa}{\pa x_i}, \frac{\pa}{\pa x_j}) = g_{ij}. $

\bel 
\label{geodesic-coordinates}
There exists a system of local coordinates $(x_1, x_2)$ centered at $p$ satisfying that 
\\
{\rm (1)} 
$g_{11} = 1, g_{12} = g_{21} = 0$ so that $g$ is of form $(d x_1)^2 + g_{22} (dx_2)^2$ and moreover 
$g_{22}$ satisfies the conditions $g_{22}(0, x_2) = 1$ and $\dfrac{\pa g_{22}}{\pa x_1}(0, x_2) = 0$. 
\\
{\rm (2)} $\dfrac{\pa g_{ij}}{\pa x_k}(p) = 0$ for any $i, j, k = 1, 2$, and all connection coefficients 
(Christoffel symbols) $\Gamma_{ij}^k$ vanish at $p$. 
\\
{\rm (3)}
$v = \frac{\pa}{\pa x_1}\vert_p$ and the curve $x_1(t) = t, x_2(t) = c$, $c$ being a constant,  gives 
the Riemannian geodesic on $M$. 
If $c = 0$, it is the Riemannian geodesic starting from $p$ with the initial velocity vector $v$. 
\\
{\rm (4)}
Let $\theta$ be the angle function with the base section $\frac{\pa}{\pa x_1}$. 
Then the generating vector field $V$ of the geodesic flow $V$ on $UM$ satisfies $\langle d\theta, V\rangle(x_1, x_2, 0) = 0$. 
\enl

\noindent
{\it Proof of Lemma \ref{geodesic-coordinates}: }
We take the geodesic parallel coordinates around $p$ on $M$. See, for instance, \cite{Kreyszig} 
on its existence. Then we have (1). The assertion (2) follows from (1) and 
$\Gamma_{ij}^k = \frac{1}{2}g^{k\ell}\left( \frac{\pa g_{\ell i}}{\pa x_j} + 
\frac{\pa g_{\ell j}}{\pa x_i} - \frac{\pa g_{i j}}{\pa x_\ell}\right)$ using Einstein convention. 
The assertion (3) follows from (1) by rotating the coordinates $(x_1, x_2)$ linearly if necessary. 
From (3), we see the $\frac{\pa}{pa \theta}$ of $V$ vanishes along the each geodesic $x_1(t) = t, x_2(t) = c$. 
Therefore we have $\langle d\theta, V\rangle(x_1, x_2, 0) = 0$, and we have (4). 
\QED

\bel
\label{orthonormal-frame}
There exists a local orthonormal frame $v_1, v_2$ 
of $TM$ on a neighborhood of $p$ such that, for some geodesic parallel coordinates 
$x_1, x_2$, they are written as 
$$
v_1 = k(x_1, x_2) \dfrac{\pa}{\pa x_1} + \ell(x_1, x_2)\dfrac{\pa}{\pa x_2}, 
\quad 
v_2 = m(x_1, x_2) \dfrac{\pa}{\pa x_1} + n(x_1, x_2)\dfrac{\pa}{\pa x_2}, 
$$
with all of first order partial derivatives of $k, \ell, m, n$ vanished at $p$. 
\enl

\Proof
In general, if we set $k = \frac{1}{\sqrt{g_{11}}}, \ell = 0, m = - \frac{g_{12}}{\sqrt{g_{11}}\sqrt{g_{11}g_{12} - g_{12}^2}}$ and $n = \frac{\sqrt{g_{11}}}{\sqrt{g_{11}}\sqrt{g_{11}g_{12} - g_{12}^2}}$, then 
$v_1 = k\frac{\pa}{\pa x_1} + \ell\frac{\pa}{\pa x_1}, 
v_2 =  m \frac{\pa}{\pa x_1} + n\frac{\pa}{\pa x_2}$ form a local orthonormal frame. 
If $(x_1, x_2)$ is a system of geodesic parallel coordinates, then we see 
$k = 1, \ell = 0, m = 0$ and $\frac{1}{\sqrt{g_{22}}}$. 
Then, for the exterior derivatives, we have $dk = d\ell = d n = 0$ and 
$dk = - \frac{k}{2g_{22}} dg_{22}$. Thus, by Lemma \ref{geodesic-coordinates} (2), we have the result. 
\QED

\section{The case of general Riemannian surfaces}
\label{The case of general Riemannian surfaces}

Let us study the case with a general Riemannian surface $(M, g)$. 
Let $\Gamma : (\R, t_0) \to UM$ be any curve-germ. 
Let $(x_1, x_2)$ be a system of geodesic parallel coordinates of $M$ centered at $\pi(\Gamma(t_0))$. 
Let 
$$
v_1 = k(x_1, x_2)\frac{\pa}{\pa x_1} + \ell(x_1, x_2)\frac{\pa}{\pa x_2}, 
\quad 
v_2 = m(x_1, x_2)\frac{\pa}{\pa x_1} + n(x_1, x_2)\frac{\pa}{\pa x_2}
$$
be a local orthonormal frame for $g$ on $M$. 
For the local coordinates $x_1, x_2, \theta$ of $UM$ introduced in \S \ref{Basic constructions from Riemannian surfaces} for the base section $\dfrac{\pa}{\pa x_1}$, 
we have that 
$$
V_1 = v_1\cos\theta  +  v_2\sin\theta, 
\quad
V_2 = \frac{\pa}{\pa \theta}, 
$$
form a local sub-Riemannian orthonormal frame of $D \subset T(UM)$. 
In this case the Hamiltonian for sub-Riemannian geodesics is given by
$$
H = u_1\{ (kp_1 + \ell p_2)\cos\theta + (m p_1 + n p_2)\sin\theta \} + u_2\varphi + \frac{1}{2}c(u_1^2 + u_2^2), 
$$
where $c \in \R$, and the equation for the extremal $(x_1, x_2, \theta, p_1, p_2, \varphi)$ is written as
\[
\left\{\ 
\begin{aligned}
& \dot{x}_1 = u_1(k\cos\theta + m\sin\theta), \quad \dot{x}_2 = u_1(\ell\cos\theta + n\sin\theta), \quad \dot{\theta} = u_2
\\
& \dot{p}_1 = - u_1\{ (\dfrac{\pa k}{\pa x_1}p_1 + \dfrac{\pa \ell}{\pa x_1}p_2)\cos\theta 
+ (\dfrac{\pa m}{\pa x_1}p_1 + \dfrac{\pa n}{\pa x_1}p_2)\sin\theta \}, 
\\
& \dot{p}_2 = - u_1\{ (\dfrac{\pa k}{\pa x_2}p_1 + \ell_{x_2}p_2)\cos\theta 
+ (\dfrac{\pa m}{\pa x_2}p_1 + \dfrac{\pa n}{\pa x_2}p_2)\sin\theta\}, 
\\
& \dot{\varphi} = - u_1\{ - (k p_1 + \ell p_2)\sin\theta + (m p_1 + n p_2)\cos\theta\}
\end{aligned}
\right.
\]
with the constraint 
$$
(kp_1 + \ell p_2)\cos\theta + (m p_1 + n p_2)\sin\theta + cu_1 = 0, \quad \varphi + cu_2 = 0. 
$$

Suppose $c = 0$.  Then, by the constraint, 
$(k p_1+\ell p_2)\cos\theta + (m p_1 + n p_2)\sin\theta = 0$ and $\varphi = 0$. For any $t$ with $u_1(t) \not= 0$, 
we have $-(k p_1 + \ell p_2)\sin\theta + (m p_1 + n p_2)\cos\theta = 0$. 
Then we have 
\[
\left(
\begin{array}{cc}
\cos\theta & \sin\theta
\\
- \sin\theta & \cos\theta
\end{array}
\right)
\left(
\begin{array}{cc}
k & \ell
\\
m & n
\end{array}
\right)
\left(
\begin{array}{c}
p_1
\\
p_2
\end{array}
\right)
= 
\left(
\begin{array}{c}
0
\\
0
\end{array}
\right)
\]
Then we have $p_1 = p_2 = \varphi = 0$, which leads a contradiction. Therefore $u_1(t)$ must be $0$ 
almost everywhere. Then $x_1(t), x_2(t), p_1(t), p_2(t)$ are locally constants and therefore $(\cos\theta(t), \sin\theta(t))$ 
and $\theta(t)$ must be a constant. 
Thus we have checked directly that any non-trivial $D$-geodesic is normal in our case.

Now suppose $c \not= 0$. By replacing $-\frac{1}{c}p_1, -\frac{1}{c}p_2, -\frac{1}{c}\varphi$ by 
$p_1, p_2, \varphi$ respectively, we may set $c = -1$. 
Then we have 
$$
u_1 =  (kp_1 + \ell p_2)\cos\theta + (m p_1 + n p_2)\sin\theta, \quad u_2 = \varphi. 
$$
Thus we have a first order ordinary differential equation 
\[
\left\{ 
\begin{aligned}
\dot{x}_1 & = \  \{ (k p_1 + \ell p_2)\cos\theta + (m p_1 + n p_2)\sin\theta\}(k\cos\theta + m\sin\theta), 
\\
\dot{x}_2 & = \ \{ (k p_1 + \ell p_2)\cos\theta + (m p_1 + n p_2)\sin\theta\}(\ell\cos\theta + n\sin\theta), 
\\
\dot{\theta} \ & = \ \varphi, 
\\
\dot{p}_1 & =    -\{ (kp_1 + \ell p_2)\cos\theta + (m p_1 + n p_2)\sin\theta\}\{ (\frac{\pa k}{\pa x_1}p_1 + \frac{\pa \ell}{\pa x_1}p_2)\cos\theta + (\frac{\pa m}{\pa x_1}p_1 + \frac{\pa n}{\pa x_1}p_2)\sin\theta \}, 
\\
\dot{p}_2 & =    -\{ (kp_1 + \ell p_2)\cos\theta + (m p_1 + n p_2)\sin\theta\}\{ (\frac{\pa k}{\pa x_2}p_1 + \frac{\pa \ell}{\pa x_2}p_2)\cos\theta 
+ (\frac{\pa m}{\pa x_2}p_1 + \frac{\pa n}{\pa x_2}p_2)\sin\theta\}, 
\\
\dot{\varphi} \ & =  -\{ (k p_1 + \ell p_2)\cos\theta + (m p_1 + n p_2)\sin\theta\}\{ -(k p_1 + \ell p_2)\sin\theta + (m p_1 + n p_2)\cos\theta \}. 
\end{aligned}
\right.
\]
In our general case, we have $\ddot{\theta} = - r\sin(2\theta + \rho)$, where 
$r = r(t) = \frac{1}{2}\{ (k p_1 + \ell p_2)^2 + (m p_1 + n p_2)^2 \}$ and 
$\rho = \rho(t)$ satisfies 
$\sin\rho = -(k p_1+ \ell p_2)(m p_1 + n p_2)/r$ and 
$\cos\rho = \frac{1}{2}\{ (k p_1+ \ell p_2)^2 - (m p_1 + n p_2)^2\}/r$. 
Note that $r$ and $\rho$ depend on $t$. Therefore we observe that $\theta$ follows a generalized equation of pendulum. 

Now we return to show Theorem \ref{Main-Theorem}. 

\ber
{\rm
It is known that any normal sub-Riemannian geodesic $(x(t), \theta(t))$ is obtained as the projection of a solution 
$x(t), \theta(t), p(t), \varphi(t)$ of the Hamiltonian equation 
\[
\dot{x} = \frac{\pa \widetilde{H}}{\pa p}, \ \dot{\theta} = \frac{\pa \widetilde{H}}{\pa \varphi}, 
\quad 
\dot{p} = - \frac{\pa \widetilde{H}}{\pa x}, \ \dot{\varphi} = - \frac{\pa \widetilde{H}}{\pa \theta}, 
\]
on $T^*(UM)$ for the another Hamiltonian 
$$
\widetilde{H}(x, \theta, p, \varphi) = \frac{1}{2} \left( \langle p, V_1\rangle^2 + \langle \varphi,  V_2\rangle^2\right). 
$$
See \cite{Montgomery} Theorem 1.14. One can check that the same result is obtained 
also by analyzing the above Hamiltonian equation. 
}
\enr

\

\noindent
{\it Proof of Theorem \ref{Main-Theorem} in the general case.}

Let $\Gamma : (\R, t_0) \to (UM, \Gamma(t_0)$ be a $D$-geodesic and 
$\widetilde{\Gamma} : (\R, t_0) \to (T^*(UM), \widetilde{\Gamma}(t_0))$ 
be a corresponding extremal for some $c \not= 0$. 
Set $\Gamma(t) = (x_1(t), x_2(t), \theta(t))$ and 
$\widetilde{\Gamma}(t) = (\Gamma(t); p_1(t), p_2(t), \varphi(t))$.  
We set 
\[
\begin{array}{rcl}
A(x_1, x_2, \theta, p_1, p_2) & := & (k p_1 + \ell p_2)\cos\theta + (m p_1 + n p_2)\sin\theta, 
\\
B(x_1, x_2, \theta, p_1, p_2) & := & -(k p_1 + \ell p_2)\sin\theta + (m p_1 + n p_2)\cos\theta. 
\end{array}
\]

{\it Part I}. $\pi$-Legendre classification of $\Gamma$. \ 
If $(\dot{x}_1(t_0), \dot{x}_2(t_0)) \not= (0, 0)$, then we have the case (iii). 
Suppose $(\dot{x}_1(t_0), \dot{x}_2(t_0)) = (0, 0)$. 
Then $A(x_1(t_0), x_2(t_0), \theta(t_0), p_1(t_0), p_2(t_0)) = 0$ at $t_0$. 
Assume $\dot{\theta}(t_0) = 0$. 
Then we have $\dot{x}_1(t_0) = \dot{x}_2(t_0) = \dot{\theta}(t_0) = \dot{p_1}(t_0) = \dot{p_2}(t_0) = \dot{\varphi}(t_0) = 0$. Therefore, by the uniqueness of solution, we see $\widetilde{\Gamma}$ itself is a constant curve, 
and $\Gamma$ is also constant, then we have (i). 
Suppose $\dot{\theta}(t_0) \not= 0$. Set $a = \dot{\theta}(t_0)$. 
If $B(x_1(t_0), x_2(t_0), \theta(t_0), p_1(t_0), p_2(t_0))  = 0$ at $t_0$, 
then we have $(p_1(t_0), p_2(t_0)) = (0, 0)$. Since 
$(p_1(t), p_2(t))$ satisfies a linear homogeneous differential equation as above, we see 
$p_1(t)$ and $p_2(t)$ are identically zero. Then $\pi\circ\Gamma$ is a constant map. 
and then $\theta(t) = at$. Thus we have the case (ii). 

Suppose $B(x_1(t_0), x_2(t_0), \theta(t_0), p_1(t_0), p_2(t_0)) \not= 0$ at $t_0$. 
Now we calculate $\Delta$ as in Lemma \ref{recognition-of-cusp}. 
\[
\begin{array}{c}
\dot{x}_1 = \left(\frac{\pa A}{\pa p_1}\right)A, \quad  \dot{x}_2 = \left(\frac{\pa A}{\pa p_2}\right) A, 
\\
\ddot{x}_1 = \left(\frac{\pa A}{\pa p_1}\right)' A + \left(\frac{\pa A}{\pa p_1}\right) A', \quad 
\ \ddot{x}_2 = \left(\frac{\pa A}{\pa p_2}\right)' A + \left(\frac{\pa A}{\pa p_2}\right) A', 
\\
\dddot{x}_1 
= \left(\frac{\pa A}{\pa p_1}\right)'' A + 2\left(\frac{\pa A}{\pa p_1}\right)' A' + \left(\frac{\pa A}{\pa p_1}\right)A'', 
\quad  
\dddot{x}_2 = \left(\frac{\pa A}{\pa p_2}\right)''A + 2\left(\frac{\pa A}{\pa p_2}\right)'A' + \left(\frac{\pa A}{\pa p_2}\right)A'',
\end{array}
\]
At $t = t_0$, we have $A = 0, \dot{k} = \dot{\ell} = \dot{m} = \dot{n} = 0, \dot{p}_1 - \dot{p}_2 = 0$. 
Therefore we have 
$$
\Delta := 
\left|
\begin{array}{cc}
\ddot{x}_1 & \dddot{x}_1 
\\
\ddot{x}_2 & \dddot{x}_2
\end{array}
\right|(t_0)
= A'(t_0) 
\left|
\begin{array}{cc}
\frac{\pa A}{\pa p_1} & 2\left(\frac{\pa A}{\pa p_1}\right)' A' + \left(\frac{\pa A}{\pa p_1}\right)A''
\\
\frac{\pa A}{\pa p_2} & 2\left(\frac{\pa A}{\pa p_2}\right)' A' + \left(\frac{\pa A}{\pa p_2}\right)A''
\end{array}
\right|(t_0)
= 
2A'(t_0)^2 
\left|
\begin{array}{cc}
\frac{\pa A}{\pa p_1} & \left(\frac{\pa A}{\pa p_1}\right)' 
\\
\frac{\pa A}{\pa p_2} & \left(\frac{\pa A}{\pa p_2}\right)' 
\end{array}
\right|(t_0)
$$
We have, at $t = t_0$, 
$$
A' = \{ -(k p_1 + \ell p_2)\sin\theta + (m p_1 + n p_2)\cos\theta\}\dot{\theta} = B\dot{\theta} 
$$
and 
$$
\left|
\begin{array}{cc}
\frac{\pa A}{\pa p_1} & \left(\frac{\pa A}{\pa p_1}\right)' 
\\
\frac{\pa A}{\pa p_2} & \left(\frac{\pa A}{\pa p_2}\right)' 
\end{array}
\right|(t_0) 
= \dot{\theta}(t_0)(k n - \ell m)(t_0)
$$
Therefore we have 
$$
\Delta = 2B(t_0)^2 (kn - \ell m)(t_0)\dot{\theta}(t_0)^3 \not= 0. 
$$
Therefore $\pi\circ\Gamma$ is right-left equivalent to the cusp and this is the case (v). 

{\it Part II}. $\pi'$-Legendre classification of $\Gamma$. \ 

Next we analyse $(\Gamma, \pi')$. 
If $\pi'\circ\Gamma$ is an immersion at $t_0$, then we have (iii). 
If $\Gamma$ is a constant curve, then we have (i). If $\pi'\circ\Gamma$ is a constant curve, then 
since $\Gamma$ is an immersion, we have (ii). 

Now suppose $\pi'\circ\Gamma$ is not an immersion at $t_0$.
We take geodesic parallel coordinates around $\pi\circ\Gamma(t_0)$ and 
local frame $v_1, v_2$ as in Lemma \ref{geodesic-coordinates}. 

Set $R(x_1, x_2, \theta) := k(x_1, x_2)\cos\theta + m(x_1, x_2) \sin\theta)$ and 
$S(x_1, x_2, \theta) := \ell(x_1, x_2)\cos\theta + n(x_1, x_2)\sin\theta)$. 
Then the geodesic flow  in a neighborhood of $\Gamma(t_0)$ is written as
$$
V = R\dfrac{\pa}{\pa x_1} + S\dfrac{\pa}{\pa x_2} 
+ W\dfrac{\pa}{\pa \theta}
$$
for some function $W = W(x_1, x_2, \theta)$. 
We use the geodesic parallel coordinates around $\pi\circ\Gamma(t_0)$. 
Then, by Lemma \ref{geodesic-coordinates}, we see $W(x_1, x_2, 0) = 0$. 

The projection $\pi'$ is locally expressed by taking a pair of independent first integrals $(F, E)$ of 
$V$ in a neighborhood of $\Gamma(t_0)$. We set 
$\pi'\circ\Gamma(t) = (F(\Gamma(t)), E(\Gamma(t))) = : (f(t), e(t))$. 
Then we have
\[
\left\{ 
\begin{array}{rcl}
\dot{f}(t) & = & \left( \dfrac{\pa F}{\pa x_1}\circ\Gamma\right)\!(t)\ \dot{x}_1(t) 
+ \left( \dfrac{\pa F}{\pa x_2}\circ\Gamma\right)\!(t)\ \dot{x}_2(t) 
+ \left( \dfrac{\pa F}{\pa \theta}\circ\Gamma\right)\!(t)\ \dot{\theta}(t), 
\\
\dot{e}(t) & = & \left( \dfrac{\pa E}{\pa x_1}\circ\Gamma\right)\!(t)\ \dot{x}_1(t) 
+ \left( \dfrac{\pa E}{\pa x_2}\circ\Gamma\right)\!(t)\ \dot{x}_2(t) 
+ \left( \dfrac{\pa E}{\pa \theta}\circ\Gamma\right)\!(t)\ \dot{\theta}(t), 
\end{array}
\right.
\]
\[
\left\{
\begin{array}{rcl}
\ddot{f}(t) & = & \left( \dfrac{\pa^2 F}{\pa x_1^2}\circ\Gamma\right)\!(t)\ \dot{x}_1(t)^2 
+ \left( \dfrac{\pa^2 F}{\pa x_2^2}\circ\Gamma\right)\!(t)\ \dot{x}_2(t)^2 
+ \left( \dfrac{\pa^2 F}{\pa \theta^2}\circ\Gamma\right)\!(t)\ \dot{\theta}(t)^2 
\\
& & 
+ 2\left( \dfrac{\pa^2 F}{\pa x_1\pa x_2}\circ\Gamma\right)\!(t)\ \dot{x}_1(t)\dot{x}_2(t) 
+ 2\left( \dfrac{\pa^2 F}{\pa x_1\pa \theta}\circ\Gamma\right)\!(t)\ \dot{x}_1(t)\dot{\theta}(t) 
+ 2\left( \dfrac{\pa^2 F}{\pa x_2\pa \theta}\circ\Gamma\right)\!(t)\ \dot{x}_2(t)\dot{\theta}(t) 
\\
& & 
+ \left( \dfrac{\pa F}{\pa x_1}\circ\Gamma\right)\!(t)\ \ddot{x}_1(t) 
+ \left( \dfrac{\pa F}{\pa x_2}\circ\Gamma\right)\!(t)\ \ddot{x}_2(t) 
+ \left( \dfrac{\pa F}{\pa \theta}\circ\Gamma\right)\!(t)\ \dot{\theta}(t), 
\\
\ddot{e}(t) & = & \left( \dfrac{\pa^2 E}{\pa x_1^2}\circ\Gamma\right)\!(t)\ \dot{x}_1(t)^2 
+ \left( \dfrac{\pa^2 E}{\pa x_2^2}\circ\Gamma\right)\!(t)\ \dot{x}_2(t)^2 
+ \left( \dfrac{\pa^2 E}{\pa \theta^2}\circ\Gamma\right)\!(t)\ \ddot{\theta}(t)^2 
\\
& & 
+ 2\left( \dfrac{\pa^2 E}{\pa x_1\pa x_2}\circ\Gamma\right)\!(t)\ \dot{x}_1(t)\dot{x}_2(t) 
+ 2\left( \dfrac{\pa^2 E}{\pa x_1\pa \theta}\circ\Gamma\right)\!(t)\ \dot{x}_1(t)\dot{\theta}(t) 
+ 2\left( \dfrac{\pa^2 K}{\pa x_2\pa \theta}\circ\Gamma\right)\!(t)\ \dot{x}_2(t)\dot{\theta}(t) 
\\
& & 
+ \left( \dfrac{\pa E}{\pa x_1}\circ\Gamma\right)\!(t)\ \ddot{x}_1(t) 
+ \left( \dfrac{\pa E}{\pa x_2}\circ\Gamma\right)\!(t)\ \ddot{x}_2(t) 
+ \left( \dfrac{\pa E}{\pa \theta}\circ\Gamma\right)\!(t)\ \ddot{\theta}(t), 
\end{array}
\right.
\]
Since $\pi'\circ\Gamma$ is not an immersion at $t_0$, we have 
$(\dot{f}(t_0), \dot{e}(t_0)) = (0, 0)$.  
Moreover we have 
$$
\dot{x}_1(t_0) = k(kp_1 + \ell p_2)\vert_{t=t_0}, \dot{x}_2(t_0) = \ell(kp_1 + \ell p_2)\vert_{t=t_0}, \dot{\theta}(t_0) = 0. 
$$
Moreover, since $\dot{p}_1(t_0) = \dot{p}_2(t_0) = 0$ and all partial derivatives of first order of 
$k, \ell, m, n$ vanish at $\\pi\circ\Gamma(t_0)$, we see 
$\ddot{x}_1(t_0) = \ddot{x}_2(t_0) = 0$. On the other hand we have $\ddot{\theta}(t_0) \not= 0$. 
Thus we have 
\[
\left\{ 
\begin{array}{rcl}
\ddot{f}(t_0) & = & \left.  \dfrac{\pa^2 F}{\pa x_1^2}k^2(kp_1+\ell p_2)^2
+ 2\dfrac{\pa^2 F}{\pa x_1\pa x_2}k\ell(kp_1+\ell p_2)^2
+ \dfrac{\pa^2 F}{\pa x_2^2}\ell^2(kp_1+\ell p_2)^2 + \dfrac{\pa F}{\pa\theta}\ddot{\theta}\ \right\vert_{t=t_0} 
\\
& = & \left. k(kp_1 + \ell p_2)^2\left( k\dfrac{\pa^2 F}{\pa x_1^2} + \ell \dfrac{\pa^2 F}{\pa x_1\pa x_2}\right) 
+ \ell(kp_1 + \ell p_2)^2\left( k\dfrac{\pa^2 F}{\pa x_1\pa x_2} + 
\ell\dfrac{\pa^2 F}{\pa x_1^2}\right) + \dfrac{\pa F}{\pa\theta}\ddot{\theta}\ \right\vert_{t=t_0}, 
\\
\ddot{e}(t_0) & = & \left.\dfrac{\pa^2 E}{\pa x_1^2}k^2(kp_1+\ell p_2)^2
+ 2\dfrac{\pa^2 E}{\pa x_1\pa x_2}k\ell(kp_1+\ell p_2)^2
+ \dfrac{\pa^2 E}{\pa x_2^2}\ell^2(kp_1+\ell p_2)^2 
+ \dfrac{\pa E}{\pa\theta}\ddot{\theta}\ \right\vert_{t=t_0} 
\\
& = & \left. k(kp_1 + \ell p_2)^2\left( k\dfrac{\pa^2 E}{\pa x_1^2} + \ell \dfrac{\pa^2 E}{\pa x_1\pa x_2}\right) 
+ \ell(kp_1 + \ell p_2)^2\left( k\dfrac{\pa^2 E}{\pa x_1\pa x_2} + 
\ell\dfrac{\pa^2 E}{\pa x_1^2}\right) + \dfrac{\pa E}{\pa\theta}\ddot{\theta}\ \right\vert_{t=t_0}, 
\end{array}
\right.
\]
Because $F$ and $E$ are first integrals of $V$, as functions on $x_1, x_2, \theta$, 
\[
\begin{aligned}
R\dfrac{\pa F}{\pa x_1} + S\dfrac{\pa F}{\pa x_2} + W\dfrac{\pa F}{\pa \theta} = 0, 
\quad 
R\dfrac{\pa E}{\pa x_1} + S\dfrac{\pa H}{\pa x_2} + W\dfrac{\pa H}{\pa \theta} = 0. 
\end{aligned}
\]
By taking the differentials by $x_1$ and $x_1$ of the right hand sides of the above equations respectively, we have 
\[
\left\{
\begin{aligned}
\dfrac{\pa R}{\pa x_1}\dfrac{\pa F}{\pa x_1} + \dfrac{\pa S}{\pa x_1}\dfrac{\pa F}{\pa x_2} 
+ \dfrac{\pa W}{\pa x_1}\dfrac{\pa F}{\pa \theta} + 
R\dfrac{\pa^2 F}{\pa x_1^2} + S\dfrac{\pa F}{\pa x_1\pa x_2} + W\dfrac{\pa F}{\pa x_1\pa \theta}= 0, 
\\
\dfrac{\pa R}{\pa x_2}\dfrac{\pa F}{\pa x_1} + \dfrac{\pa S}{\pa x_2}\dfrac{\pa F}{\pa x_2} 
+ \dfrac{\pa W}{\pa x_2}\dfrac{\pa F}{\pa \theta} 
+ 
R\dfrac{\pa^2 F}{\pa x_1\pa x_2} + S\dfrac{\pa F}{\pa x_2^2} + W\dfrac{\pa F}{\pa x_2\pa \theta}
= 0, 
\\
\dfrac{\pa R}{\pa x_1}\dfrac{\pa E}{\pa x_1} + \dfrac{\pa S}{\pa x_1}\dfrac{\pa E}{\pa x_2} 
+ \dfrac{\pa W}{\pa x_1}\dfrac{\pa E}{\pa \theta} 
+ R\dfrac{\pa^2 E}{\pa x_1^2} + S\dfrac{\pa E}{\pa x_1\pa x_2} + W\dfrac{\pa E}{\pa x_1\pa \theta}
= 0, 
\\
\dfrac{\pa R}{\pa x_2}\dfrac{\pa E}{\pa x_1} + \dfrac{\pa S}{\pa x_2}\dfrac{\pa E}{\pa x_2} 
+ \dfrac{\pa W}{\pa x_2}\dfrac{\pa E}{\pa \theta} 
+ 
R\dfrac{\pa^2 E}{\pa x_1\pa x_2} + S\dfrac{\pa E}{\pa x_2^2} + W\dfrac{\pa E}{\pa x_2\pa \theta} = 0, 
\end{aligned} 
\right.
\]
At the point $\Gamma(t_0)$, we have all of 
$\dfrac{\pa R}{\pa x_1}, \dfrac{\pa R}{\pa x_2}, \dfrac{\pa S}{\pa x_1}, \dfrac{\pa S}{\pa x_2}, 
W, \dfrac{\pa W}{\pa x_1}, \dfrac{\pa W}{\pa x_2}$ vanish. 
Moreover we have $R(\Gamma(t_0)) = k(\pi\circ\Gamma(t_0)), S(\Gamma(t_0) = \ell(\pi\circ\Gamma(t_0))$. 
Therefore we obtain 
\[
\begin{aligned}
\left. k\dfrac{\pa^2 F}{\pa x_1^2} + 
\ell \dfrac{\pa^2 F}{\pa x_1\pa x_2}\right\vert_{t = t_0} = 0, 
\quad 
\left. k\dfrac{\pa^2 F}{\pa x_1\pa x_2} + 
\ell\dfrac{\pa^2 F}{\pa x_1^2} \right\vert_{t = t_0} = 0, 
\\
\left. k\dfrac{\pa^2 E}{\pa x_1^2} + 
\ell \dfrac{\pa^2 E}{\pa x_1\pa x_2}\right\vert_{t = t_0} = 0, 
\quad 
\left. k\dfrac{\pa^2 E}{\pa x_1\pa x_2} + 
\ell\dfrac{\pa^2 E}{\pa x_1^2}\right\vert_{t = t_0} = 0. 
\end{aligned}
\]
Therefore we have 
$$
\ddot{f}(t_0) = \dfrac{\pa F}{\pa \theta}(t_0)\ddot{\theta}(t_0), \quad 
\ddot{h}(t_0) = \dfrac{\pa E}{\pa \theta}(t_0)\ddot{\theta}(t_0). 
$$
Since $\dfrac{\pa}{\pa \theta}$ does not belongs to the kernel of the differential of $(F. E) : (UM, \Gamma(t_0)) 
\to \R^2$ at $\Gamma(t_0)$, we have $(\dfrac{\pa F}{\pa \theta}(t_0), \dfrac{\pa E}{\pa \theta}(t_0)) 
\not= (0, 0)$. Since $\ddot{\theta}(t_0) \not= 0$, we have 
$(\ddot{f}(t_0), \ddot{e}(t_0)) \not= (0, 0)$. Therefore, by Lemma \ref{recognition-of-cusp-2}, 
we have $\pi'\circ\Gamma$ is a cusp. 

The last statement on the combination of singularities of $\pi\circ\Gamma$ and $\pi'\circ\Gamma$ 
follows by remarking that the combinations ((iii), (iii)), ((iv), (iv)) are never occur because 
$\Ker(d\pi) \cap \Ker(d\pi') = V \cap K = \{ 0\}$. 
This completes the proof of Theorem \ref{Main-Theorem}.

\section{Appendix: A naive motivation} 
\label{Appendix: A naive motivation}

In the winter snow season, you will observe many cusp-shaped traces of vehicles on many roads and parking lots usually. 
Naturally it can be supposed that we control vehicles in a (nearly) optimal way, when we drive and park. 
Therefore the cuspidal shape of such snow-traces may be regarded as an appearance of generic singularities for 
solutions to some problem of optimal control theory. For instance: 

\

Problem. Suppose your car is located on a parking place. You are asked to move your car to the very next (right) place. 
How do you drive and move your car ? 

\begin{center}
\includegraphics[width=9truecm, height=3truecm, clip, 
]{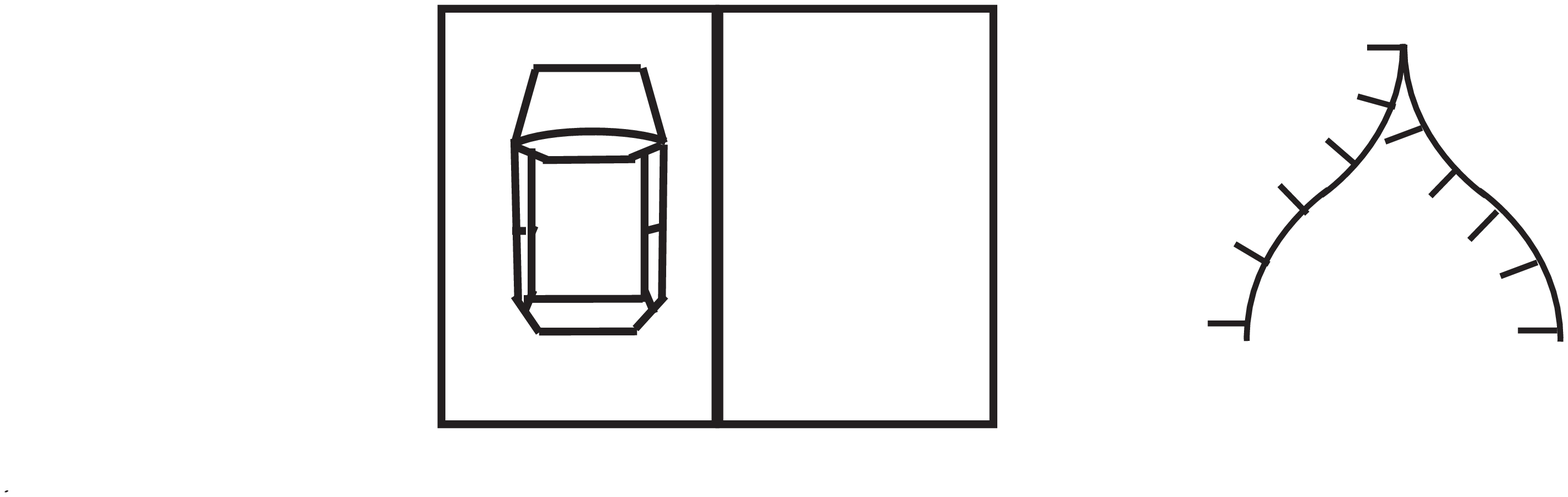} 
\end{center}

Maybe you will go forward to a right direction a little and then go back to the proper parking place. 
Then the trace of your drive wheel will form a cusp-shaped curve. 
The front direction of the wheel or its left-side normal is determined anytime, so the trace can be regarded 
as a kind of so-called a \lq\lq front" or a \lq\lq frontal" (\cite{Ishikawa18}). The short lines indicate the left side 
directions of the driver, which form a normal field to the trace.

The phenomena of the appearance of singularities do not depend on the flatness of the field and 
you will observe the singularities also on slopes and non-flat parking lots everywhere. 
The singularities can be understood as Legendre singularities of sub-Riemannian geodesics for 
general Riemannian surfaces which we have discussed in the present paper.

\

\begin{flushleft}

(G. Ishikawa)
Department of Mathematics, Faculty of Science, Hokkaido University, Japan. 
\\
E-mail: ishikawa@math.sci.hokudai.ac.jp

\

(Y. Kitagawa)
Oita National College of Technology, Oita 870-0152, Japan. 
\\
E-mail: kitagawa@oita-ct.ac.jp

\end{flushleft}

\end{document}